\documentclass[11pt,a4paper,reqno]{amsart}

\usepackage{amsfonts,amsthm,amscd,epsfig,amsmath,amssymb,enumerate}

\author {A. Faggionato}
\address{Dipartimento di Matematica ``G. Castelnuovo", Universit\`a   ``La
  Sapienza'', P.le Aldo Moro  2, 00185  Roma, Italy.}
\email{faggiona@mat.uniroma1.it}

\newtheorem{theo}{Theorem}

\newtheorem{prop}{Proposition}
\newtheorem{lemma}{Lemma}

\newtheorem{rem}{Remark}

\newcommand{\NN}{{\mathbb N}}
\newcommand{\RR}{{\mathbb R}}
\newcommand{\ZZ}{{\mathbb Z}}

\newcommand{\QQ}{{\mathbb Q}}
\newcommand{\PP}{{\mathbb P}}
\newcommand{\EE}{{\mathbb E}}

\newcommand{\Aa}{{\mathcal A}}

\newcommand{\Bb}{{\mathcal B}}
\newcommand{\Gg}{{\mathcal G}}

\newcommand{\Ee}{{\mathcal E}}
\newcommand{\Ff}{{\mathcal F}}

\newcommand{\Cc}{{\mathcal C}}

\newcommand{\Ll}{{\mathcal L}}

\newcommand{\II}{{\mathbb I}}
\let\a=\alpha \let\b=\beta   \let\d=\delta  \let\e=\varepsilon
 \let\g=\gamma \let\h=\eta

\let\o=\omega
\let\p=\pi  

\let\s=\sigma  
\let\t=\tau

\let\O=\Omega  \let\P=\Pi

\title{Bulk diffusion of 1D exclusion process with bond disorder.}

\begin{document}

\maketitle

\begin{abstract}
Given a doubly infinite  sequence of positive numbers $\{c_k\,:\,
k\in \ZZ \}$ such that $\{c_k^{-1}\,:\, k\in \ZZ\}$ satisfies a  LLN
with limit $\a\in (0,\infty]$, we consider the nearest--neighbor
  simple exclusion process on $\ZZ$ where $c_k$ is the probability
  rate of  jumps between   $k$ and   $k+1$. If $\a=\infty$ we
  require an additional minor technical condition.
   By extending a method developed in
  \cite{N}  we show that the diffusively rescaled
    process has hydrodynamic behavior described by the heat equation
    with diffusion constant $1/\a$. In particular, the process has
    diffusive behavior for $\a<\infty $ and subdiffusive behavior for
    $\a=\infty$.\\


\end{abstract}

\section{Introduction}

We consider a particle  system  on $\ZZ$ with site exclusion
interaction  performing a   stochastic dynamics with  Markov generator
\begin{equation}\label{generare}
\Ll f(\eta)=\sum _{k\in \ZZ} c_k \left( f\left(\eta^{k,k+1}\right)-
f(\eta)\right),
\end{equation}
where $f$ is a cylinder function on the state space $\{0,1\}^{\ZZ}$
and,
 given $\eta\in
\{0,1\}^\ZZ$, $\eta ^{k,k+1}$ is defined as
$$
\eta ^{k,k+1}_x=
\begin{cases}
\eta _{k+1} & \text{ if }x=k,\\
\eta _k & \text{ if } x=k+1,\\
\eta _x & \text{ otherwise}.
\end{cases}
$$
The  family $\{c_k\}_{k\in \ZZ}$ is thought of as  the environment
of the above exclusion process. We assume that $c_k>0$ for all $k\in
\ZZ$ and that for a suitable constant $\a \in (0,\infty]$
\begin{equation}\label{asso}
\lim _{k\uparrow\infty}\frac{ S (\lfloor y k \rfloor) - S(\lfloor x
k\rfloor) }{(y-x)k}
=\a\,\qquad \forall x<y,
\end{equation}
where $\lfloor \cdot \rfloor $ denotes the integer part and the
function $S: \ZZ\to\RR$ is defined as
\begin{equation}\label{sos}
 S(k)=
\begin{cases}
\sum _{j=0}^{k-1}  \frac{1}{c_j} & \text{ if } k\geq 1 ,\\
0 & \text{ if } k=0,\\
-\sum_{j=1}^{-k} \frac{1}{c_{-j}} &\text{ if } k<0.
\end{cases}
\end{equation}
 Trivially (\ref{asso})
implies that
\begin{equation}\label{condi}
\lim _{k\uparrow\infty}\frac{1}{k}\sum _{j=0}^{k-1} \frac{1}{c_j}=\a
\,,\qquad
\lim _{k\uparrow\infty}\frac{1}{k}\sum _{j=1}^{k}
\frac{1}{c_{-j}}=\a.
\end{equation}
If $\a\in (0,\infty)$,  then (\ref{asso}) and (\ref{condi}) are
equivalent. If $\a =\infty$, then (\ref{condi}) does not imply
(\ref{asso}) (see Appendix  \ref{aaa}).

Our main results concern the hydrodynamic behavior of  the above
exclusion process. In what follows, given a probability measure $\mu$
on $\{0,1\}^\ZZ$, we denote by $\PP_\mu$ the law of  the exclusion
process with generator (\ref{generare}) and initial distribution
$\mu$.
\begin{theo}\label{Hyd}
 Suppose that $\{c_k\}_{k\in \ZZ}$  satisfies condition
(\ref{condi}) with $\a\in (0,\infty)$. Let  $\rho_0:\RR\to
[0,\infty)$ be a bounded Borel function  and let
$\left\{\mu_n\right\}_{n\geq 0}$ be a family of  probability
measures on $\{0,1\}^\ZZ$ such that, for all $\varphi\in C_0
(\RR)\cap L^1(\RR)$ and $\d>0$,
$$\lim _{n\uparrow \infty}
\mu_n \left( \left| \frac{1}{n} \sum _{j\in \ZZ} \varphi
\left(\frac{j}{n}\right) \eta _j - \int _\RR \varphi(x) \rho_0 (x)
dx \right|>\delta \right)=0.
$$
Then, for all $t>0$, $\varphi\in C_c(\RR)$ and $\d>0$,
\begin{equation}
\lim _{n\uparrow \infty} \PP_{\mu_n}
 \left( \left| \frac{1}{n} \sum _{j\in \ZZ} \varphi
\left(\frac{j}{n}\right) \eta _j(n^2 t) - \int _\RR \varphi(x) \rho
(x,t)  dx \right|>\delta \right)=0,
\end{equation}
where
 $\rho:\RR \times [0,\infty)\to \RR$ solves
the heat equation
$$
\partial _t \rho= (1/\a)\partial^2_x \rho
$$
with boundary condition $\rho _0$ at time $t=0$.
\end{theo}
We remark that the above function $\rho$ can be represented as
\begin{equation}
\rho (x,t)=\int _\RR p(t,x-y) \rho_0 (y) dy ,
\end {equation}
where  $p(t,x)$ is the density of a Gaussian variable $\mathcal{N}
(0, 2t /\a )$.

Note that  if   $\{c_k\}_{k\in \ZZ}$  are i.i.d. positive random variables such that
$\EE(1/c_k)<\infty$, then the above theorem holds for almost all realizations of  $\{c_k\}_{k\in \ZZ}$
with $\a= \EE(1/c_k)$.

\medskip

In order to discuss the subdiffusive behavior of the system it is
convenient to introduce the following notation: we say that
condition (H) is fulfilled  if for all  $x\in \QQ$, $ a\not = 0 $, $
\e >0 $ there exists a sequence  of integer numbers  $\{b_n\}_{n\geq
1}$  such that  $a b_n \geq 0 $ and
\begin{align}
& \lim _{n\uparrow \infty} \frac{b_n+|a|/a}{n}\, \frac{
S\left( \lfloor a n \rfloor +\lfloor x n \rfloor \right)
-
S\left( \lfloor x n \rfloor \right)
}{n} =\infty ,\label{ca1}\\
& \varlimsup _{n\uparrow \infty} \frac{
S\left(  b_n   +\lfloor x n \rfloor \right)
-
S\left( \lfloor x n \rfloor \right)
}{
S\left( \lfloor a n \rfloor +\lfloor x n \rfloor \right)
-
S\left( \lfloor x n \rfloor \right)
} \leq \e \label{ca2}.
\end{align}

Note that condition (H) is satisfied if (\ref{asso}) is true with
$\a=\infty$ and the following holds for all $x\in \QQ$, $a\not =0$:
\begin{equation}
\lim _{\g\downarrow 0} \varlimsup _{n\uparrow\infty}
\frac{ S\left(  \lfloor \g a  n \rfloor +\lfloor x n
\rfloor \right)
-
S\left( \lfloor x n \rfloor \right)
}{
S\left( \lfloor a n \rfloor +\lfloor x n \rfloor \right)
-
S\left( \lfloor x n  \rfloor \right) }=0 .
\end{equation}
Moreover, note that condition (H) implies   (\ref{asso})
with $\a=\infty$.

\begin{theo}\label{lento}
Suppose that $\{c_k\}_{k\in \ZZ}$  satisfies  condition (H). Let
$\rho_0:\RR\to [0,\infty)$ be in  $L^1_{\text{loc}}(\RR)$
and let $\left\{\mu_n\right\}_{n\geq 0}$ be a family of probability
measures on $\{0,1\}^\ZZ$ such that, for all $\varphi\in C_c (\RR)$
and $\d>0$,
\begin{equation}\label{ricamo}
\lim _{n\uparrow \infty} \mu_n \left( \left| \frac{1}{n} \sum _{j\in
\ZZ} \varphi \left(\frac{j}{n}\right) \eta _j - \int _\RR \varphi(x)
\rho_0 (x) dx \right|>\delta \right)=0. \end{equation}
 Then, for all
$t>0$, $\varphi\in C_c(\RR)$ and $\d>0$,
\begin{equation}
\lim _{n\uparrow \infty} \PP_{\mu_n}
 \left( \left| \frac{1}{n} \sum _{j\in \ZZ} \varphi
\left(\frac{j}{n}\right) \eta _j(n^2 t) - \int _\RR \varphi(x)
\rho_0 (x)  dx \right|>\delta \right)=0.
\end{equation}
\end{theo}
We point out that  condition (H) enters only in the proof of
Proposition \ref{supervacanza}, which can be obtained under  weaker
conditions.  For example,   as discussed in Remark  \ref{jolly}, it
is enough that  there exists  an increasing sequence of positive
integers $\{n_k\}_{k\geq 1}$ such that
\begin{equation}\label{tordo}
\lim _{k\uparrow \infty} \frac{n_{k+1}-n_k}{n_k}=0
\end{equation}
and such that for all  $x\in \QQ$, $ a\not = 0 $, $ \e >0 $  one can
define  a sequence  of integer numbers   $\{b_k\}_{k\geq 1}$  such
that $a b_k \geq 0 $ and
\begin{align}
& \lim _{k\uparrow \infty} \frac{b_k+|a|/a}{n_k}\, \frac{
S\left( \lfloor a n_k \rfloor +\lfloor x n_k \rfloor \right)
-
S\left( \lfloor x n_k \rfloor \right)
}{n_k} =\infty ,\label{cca1}\\
& \varlimsup _{k\uparrow \infty} \frac{
S\left(  b_k   +\lfloor x n_k \rfloor \right)
-
S\left( \lfloor x n_k \rfloor \right)
}{
S\left( \lfloor a n_k \rfloor +\lfloor x n_k \rfloor \right)
-
S\left( \lfloor x n_k \rfloor \right)
} \leq \e \label{cca2}.
\end{align}
In particular, whenever  the above condition is fulfilled the particle system has subdiffusive behavior (note
 that in this case condition (\ref{asso}) with $\a=\infty$ is satisfied). As example of application we prove
 in Section \ref{isa1}  the following result:
\begin{prop}\label{lentone}
Suppose that $\{c_k\}_{k\in \ZZ}$  are i.i.d. positive random
variables in the domain of attraction of a $\nu$--stable law with
$0<\nu<1$. Then the particle system is subdiffusive.
\end{prop}

The above exclusion process with bond disorder is an example  of
random barrier model (a small transition rate $c_k$ corresponds  to
a barrier between sites $k$ and $k+1$) and it has been used by
physicists to model transport of charge carriers in one dimensional
disordered media (see for example \cite{ABSO}, \cite{DS}).  From a
physical viewpoint,   (\ref{asso})   is the natural condition in
order to observe a diffusive behavior possibly with zero diffusion
constant: the diffusively rescaled process can be associated to a 1D
resistor network with $\ZZ/n$ as vertex set such that the bond
$[j/n, (j+1)/n]$ has resistance $1/( n c_j)$. Then the total
resistance of the filament  $(x,y]$ is given by
$$
 \sum _{j=\lfloor x n\rfloor+1 }^{\lfloor y n \rfloor } \frac{1}{n
c_j}
$$
and due to (\ref{asso}) it converges to $\a(y-x)$ as $n\uparrow
\infty$. Therefore, assumption (\ref{asso}) means that the linear
filament has uniform (macroscopic) resistance per unit length equal
to $\a$.
 In particular, it is
natural  to have a non trivial diffusive behavior if $\a <\infty$
and a null diffusive behavior if $\a =\infty$ (condition (H) is a
more technical condition, used only in the proof of Proposition
\ref{supervacanza}).

 Due to the above observation the conditions required in \cite{N}[Theorem 3] appear artificial.
There,  K. Nagy  proves the same result as in  Theorem \ref{Hyd}
above for almost all realization of a i.i.d. random environment
$\{c_k\}_{k\in
 \ZZ}$ by requiring that   $ E( c_k ^{-4} )<\infty$ and that $c_k \leq C<\infty$ a.s. The strategy
 followed by K. Nagy consists in showing that, for what concerns bulk diffusion,
 one can ignore the site exclusion constraint in the diffusive limit. In these notes
we   show how to improve this method by using a classing result of
C. Stone \cite{S} allowing to represent the   random walk on $\ZZ$
having $c_k$ as  probability rate of  jumps between $k,k+1$  as a
space--time change of a 1D Brownian motion (see also \cite{KK}).

We observe that by techniques which are  standard  for non gradient
systems one can prove the hydrodynamic limit for the
nearest--neighbor exclusion process on $\ZZ^d$ with bond disorder,
where $c_{x,y}$ is the probability rate  for a jump between adjacent
sites $x,y$ and $\{ c_{x,y}\,:\, |x-y|=1\} $ is a family of i.i.d.
random variables such that $0<C\leq c_{x,y} \leq C^{-1}$ a.s. The
hydrodynamic limit holds for almost any realization of the disorder
$\{ c_{x,y}\,:\, |x-y|=1\}$ and is independent from the disorder.
See for example \cite{FM} (here the canonical expectation of the
gradient density field is zero, thus simplifying drastically the
treatment in \cite{FM} and allowing to get easily a proof for any
dimension).

 We point out that the results of Section \ref{extnagy}
are valid in all dimensions, while Stone's method (treated in
Section \ref{secpietra}) works only in dimension one.
  In particular,
 the method described here allows to prove the hydrodynamic limit of
the exclusion process with bond disorder in any dimension $d$ when
having   results on the single random walk  similar to the ones
described in  Section  \ref{secpietra}. For a more detailful
discussion see Appendix \ref{alto}.

 The   hydrodynamic behaviour of one--dimensional stochastic
  processes with disorder have been studied in  several papers (e.g. \cite{Fritz}, \cite{Landim}).
For a  discussion on the hydrodynamic limit of lattice gases
 with  site disorder see \cite{FM} and references therein.

The paper is structured as follows.
   In Section \ref{lego} we show
that the dynamics of the above exclusion process is well defined and
recall its graphical representation.
 In Section \ref{extnagy} we recall and extend the method developed in
 \cite{N}[Section 4]. In Section \ref{secpietra} we study   the symmetric random walk on $\ZZ$
  with rates
$\{c_k\}_{k\in \ZZ}$ using Stone's representation.  In Section
\ref{fine} we give the proof of Theorem \ref{Hyd} and Theorem
\ref{lento}, while in Section \ref{isa1} we give the proof of
Proposition \ref{lentone}.  Finally in Appendix \ref{aaa} we show
that condition (\ref{asso}) is not equivalent to condition
(\ref{condi}) if  $\a=\infty$ and in Appendix \ref{alto} we   show
some extensions of our results to higher dimension.


\section{Graphical representation of the exclusion
process}\label{lego}

 By means  of the  graphical representation of exclusion processes
\cite{D} \cite{L2}, we prove in this section that  the dynamics of
the exclusion process  with generator (\ref{generare}) is well
defined  since  (\ref{condi}) holds with $\a\in(0,\infty]$. The
graphical representation explained below  will be used also  in
Section \ref{extnagy}.

Let  $N_k (\cdot) $, $k\in \ZZ$,  be a family of independent Poisson
processes  defined on some probability space $(\O,\Ff , P)$ such
that $E(N_k(t))=c_k t$. Given $t>0$ we define $\Gg_t$ as the random
graph with vertex set $\ZZ$ and edges $\{k,k+1\}$ such that
$N_k(t)\geq 1$.
\begin{lemma}
For almost all $\o$,   the graph $\Gg_t (\o)$ has only finite
connected components for all $t>0$.
\end{lemma}
\begin{proof}
We claim that $ P(N_x(t) \geq 1\; \forall x \geq k) =0$ for all
$k\in \ZZ$. In fact, since $1-z\leq e^{-z}$ for all $z\geq 0$,
$$ P(N_x(t) \geq 1\; \forall x \geq k)= \lim _{N\uparrow \infty}\P
_{x=k} ^N \left( 1- e^{-c_x t} \right)\leq \lim _{N\uparrow \infty}
\exp\left\{ -\sum _{x=k} ^N e^{-c_x t}
 \right\}
$$
and the sum in the last member  goes to $\infty $ as $N\uparrow
\infty$ since it cannot hold $\lim _{x\uparrow\infty} c_x =\infty$
due to (\ref{condi}).

Similarly one can prove that  $ P(N_x(t) \geq 1\; \forall x \leq  k)
=0$ for all $k\in \ZZ$. In particular, almost surely  for all $t\in
\NN$ the set $\{x\,:\, N_x (t)=0 \}$ is unbounded from below and
from above, thus implying that $\Gg_t$ has only finite connected
components for all $t\in \NN$. To conclude the proof it is enough to
observe that $\Gg_s\subset \Gg_t$ for $s\leq t$.
\end{proof}

Let $\Aa\in \Ff $, with $P(A)=1$,  be a set of configurations $\o$
such that $\Gg_t(\o)$ has only finite connected components  for all
$t>0$  and $N_k(\cdot)$ has only jumps of value  $1$ for all $k\in
\ZZ$. Let $\o\in \Aa $. Then, given an initial configuration
$\h(0)$, the configuration  $\h (t)= \h (t)[\o]$ at time $t$ is
defined  as follows:

   Let $\Cc$ be any connected
component of $ \Gg_t(\o)$ and let
$$
\{s_1<s_2<\cdots <s_r\}= \{s\,:\, N_k(s)= N_k(s-)+1 ,\; \{k,k+1\}\in
 \Cc,\; 0<s\leq t\}.
 $$
 Start with $\h(0)$. At time $s_1$ switch the values between $\h_k$
 and $\h_{k+1}$ if $N_k (s_1)=N_k(s_1-)+1$ and $\{k,k+1\}\in \Cc$.
 Repeat the same operation orderly for times $s_2,s_3, \dots , s_r$.
Then the  resulting configuration coincides with $\h(t)$ on $\Cc$.


\section{Site exclusion constraint}\label{extnagy}

 Following the main ideas of \cite{N}[Section 4], we prove in this section that
 the site exclusion constraint becomes negligible when considering
 the bulk diffusion of the particle system, i.e. from a hydrodynamic viewpoint  the system behaves
 as a family of independent continuous--time random walks on $\ZZ$ with Markov generator
   $H:\RR^\ZZ\to \RR^\ZZ$  defined as
\begin{equation}\label{genRW}
(Hf)_k = c_k\left( f_{k+1}-f_{k} \right)+ c_{k-1} \left( f_{k-1}-
f_{k} \right).
\end{equation}
Note that  the random walk on $\ZZ$ with Markov generator $H$ is
reversible since the transition rates are bond dependent. In
particular, $p(t,j,k)=p(t,k,j)$ where $p(t,x,y)$ denotes the
probability that the random walk starting at $x$ is in $y$ at time
$t$.

Since
$$ d\h_k(t)= (\h_{k+1}-\h_k)(t-) dN_k (t) +  (\h_{k-1}-\h_k)(t-) dN_{k-1}
(t),
$$
we can write
\begin{equation}\label{derivo}
 d\eta (t)= H \eta (t) dt  + d M(t)
 \end{equation}
where
$$
dM_k(t)=(\h_{k+1}-\h_k)(t-) d A_k (t) + (\h_{k-1}-\h_k)(t-) d
A_{k-1} (t)
$$
and
$$ A_x (t) = N_x (t) - c_x t.
$$
Note that
 $M_k(\cdot)$ has  trajectories of bounded
variation on finite intervals a.s.

 Formally, (\ref{derivo}) implies that
\begin{equation}\label{formale}
\h (t) = T(t) \h (0) + \int_0 ^t T(t-s) dM(s)
\end{equation}
where $T(t) = e^{tH}$, i.e.
\begin{equation}\label{formale1}
\h _k (t) = \sum_{j\in \ZZ} p(t,k,j) \h _j (0) + \sum _{j\in \ZZ}
\int _0 ^t p(t-s,k,j) dM_j (s).
\end{equation}
Due to the graphical construction of the dynamics discussed in
Section \ref{lego}, if $\sum _{x\in \ZZ}\h_x(0)<\infty$, then for
all but a finite family of indexes $j$
 $ dM_j(s)=0$ for all $0\leq s\leq t$ and in particular  the last series
in (\ref{formale1}) reduces to a finite sum and is  meaningful. In
this case, one can check that
 (\ref{formale1}) holds a.s. by direct computation
 using that
$$
\frac{d}{dt} p(t, k,j)= \left( H p(t, \cdot, j)\right)_k.
$$
The following result shows that the site exclusion constraint is
negligible in the diffusive limit from a  hydrodynamic viewpoint:
\begin{prop}\label{gaetano} Given  $\d>0$, $t>0$, $\varphi \in C_c
(\RR)$ and given a sequence of probability measures $\mu_n$ on
$\{0,1\}^\ZZ$,
\begin{equation}\label{limitone}
\lim _{n\uparrow\infty} \PP _{\mu_n} \left( \left| \frac{1}{n} \sum
_{k\in \ZZ} \varphi\left(\frac{k}{n}\right) \h _k (tn^2)-
\frac{1}{n} \sum _{k\in \ZZ} \varphi\left(\frac{k}{n}\right)\sum
_{j\in \ZZ} p(t n^2, k,j) \h _j (0) \right|>\d \right)=0.
\end{equation}
\end{prop}
\begin{proof}
 Let the support of $\varphi$ be included in $[-L,L]$ and fix $\e>0$.
Given $x\in \ZZ$ and $t>0$ define $\Cc _x (t)$ as the connected
component of $\Gg _t$ containing $x$.  Then for each positive
integer $n $ we can choose $b_n$ large enough such that $ P(\Aa _n
^c ) <\epsilon$ where $\Aa_n$ is the subset of configurations $\o$
satisfying the
following conditions:\\
\begin{align}
& \cup _{x\in [-Ln,Ln] } C_x (tn ^2) \subset [-b_n, b_n],\label{roma1}\\
& (2L+1) \|\varphi \|_\infty \sup _{k\in [-Ln, Ln] } \sum _{j\,: \,
|j|> b_n } p (tn^2,k, j ) \leq \d/2 .\label{roma2}
\end{align}
Given $\h (0)$ and $n$, we define $\h^{(n)} (0) \in \{0,1\}^\ZZ$ as
$$
\h^{(n)} _k ( 0) = \h _k (0)\II_{|k| \leq b_n}
$$
and write  $\h ^{(n)}(s)$ for  the configuration at time $s$
obtained by the graphical construction when starting from
$\h^{(n)}(0)$ at time $0$.

 Due to the graphical construction of the dynamics and
condition (\ref{roma1}), if $\o\in \Aa_n$ then
$$
\h _k  (t n^2) = \h _k ^{(n)} (t n^2),\qquad \forall k\in [-Ln, Ln],
$$
thus implying
$$
 \frac{1}{n} \sum
_{k\in \ZZ} \varphi\left(\frac{k}{n}\right) \h _k (tn^2)=
 \frac{1}{n} \sum
_{k\in \ZZ} \varphi\left(\frac{k}{n}\right) \h _k ^{(n)} (tn^2).
$$

 Moreover, due to
(\ref{roma2}), if $\o\in \Aa_n$ then
$$
 \frac{1}{n} \sum
_{k\in \ZZ} \left|\varphi\left(\frac{k}{n}\right)\right|\sum _{j\in
\ZZ} p(tn^2, k,j) |\h _j(0)- \h _j ^{(n)} (0)|\leq \d/2. $$
Therefore  the l.h.s. of (\ref{limitone}) with fixed $n$ can be
bounded by
$$ \PP_{\mu_n} \left(  | Z_n |>\d /2 \right)+ P ( \Aa ^c_n ) \leq
4 \EE_{\mu_n} \left( Z_n ^2 \right)/\d^2 +\e
$$
 where
$$
Z_n =
 \frac{1}{n} \sum
_{k\in \ZZ} \varphi\left(\frac{k}{n}\right) \h ^{(n)} _k (tn^2)-
\frac{1}{n} \sum _{k\in \ZZ} \varphi\left(\frac{k}{n}\right)\sum
_{j\in \ZZ} p(t n^2, k,j) \h _j^{(n)} (0).
$$
Since $\sum _{x\in \ZZ} \h ^{(n)} _x (0) <\infty$, setting
$$
dM^{(n)} _k(s)=(\h^{(n)} _{k+1}-\h^{(n)} _k)(s-) d A_k (s)  +
(\h^{(n)}_{k-1}-\h^{(n)}_k)(s-) d A_{k-1}(s), $$
 $(\ref{formale1}) $ implies that
 \begin{equation*}
Z_n  =  \frac{1}{n} \sum _{k\in \ZZ} \varphi\left(\frac{k}{n}\right)
\sum _{j\in \ZZ} \int _0 ^{tn^2} p(tn^2-s,k,j) dM^{(n)}_j (s).
\end{equation*}
In order to conclude the proof it is enough to apply Lemma
\ref{mezzo} to the above estimates.
\end{proof}

\begin{lemma}\label{mezzo}
For each $n\in \NN$ let $\nu _n $ be a probability measure on
$\{0,1\}^\ZZ$ such that $ \nu _n (\sum _{k\in\ZZ} \h_k <\infty) =
1$. Then
\begin{equation*}
 \lim _{n\uparrow\infty}
\EE _{\nu _n} \left(\left(
\frac{1}{n} \sum _{k\in \ZZ} \varphi\left(\frac{k}{n}\right) \sum
_{j\in \ZZ} \int _0 ^{tn^2} p(t n^2 -s,k,j) dM _j (s)
\right)^2 \right)=0.
\end{equation*}
\end{lemma}
Recall that the above series over $j$ reduces to a finite sum
whenever $\sum _{k\in \ZZ} \h_k (0) <\infty$, and therefore it is
well defined a.s.
\begin{proof}
We define $f_n$ as
\begin{multline}\label{mozart}
f_n= \frac{1}{n} \sum _{k\in \ZZ} \varphi\left(\frac{k}{n}\right)
 \sum _{j\in \ZZ} \int _0 ^{tn^2} p(tn^2 -s,k,j) dM _j (s) =\\
 \frac{1}{n} \sum _{k\in \ZZ} \varphi\left(\frac{k}{n}\right)  \sum _{j\in
\ZZ} \int _0 ^{tn^2}\left( p (tn^2-s,k, j) - p(t n^2 -s,k,
j+1)\right) (\h
  _{j+1} - \h  _j ) (s- ) d A _j (s).
\end{multline}
We remark that due to the graphical representation of the exclusion
process, $f_n$ can be thought of as a function on the probability
space $\left( \{0,1\}^\ZZ\times \O, \Bb\times \Ff, \nu_n\otimes P
\right)$, where $\Bb$ denotes the Borel $\s$--algebra of
$\{0,1\}^\ZZ$. Moreover, note that  $|f_n|\leq c(\varphi)$ due to
(\ref{formale1}).

In the following  arguments $n$ can be thought of  as fixed. Due to
our assumption on $\nu_n$, given  $\e$ with $0<\e<1$ there exists
$\ell _n\in \NN $ such that $\nu _n (A^c)\leq \e$ where
$$
A  = \left\{ \h\,:\, \h_{x} =0\text{ if } |x|\geq \ell_n \right\}.
$$
Moreover, one can find $M\in \NN$ such that $P(B ^c) \leq \e$ where
$$
B =\left\{\o\,:\, \cup _{|x|\leq \ell _n } C_x (tn^2)[\o]\subset
(-M, M) \right\}.
$$
 Then $(\nu_n \otimes P) (A\times B) \geq (1-\e)^2$. Due to the graphical
representation, one gets
$$ \II_{D} f_n = \II _{D} f_{M,n} $$
where $D=A\times B$ and
\begin{multline*}
f_{M,n}= \\
\frac{1}{n} \sum _{k\in \ZZ} \varphi\left(\frac{k}{n}\right) \sum
_{j\,: \,|j|\leq M} \int _0 ^{tn^2}\left( p (tn^2 -s,k, j) - p(t n^2
-s,k, j+1)\right) (\h _{j+1} - \h _j ) (s- ) d A _j (s).
\end{multline*}
In particular,
$$
\EE_{\nu_n} (f_n ^2)\leq c(\varphi)^2 \left(\nu _n \times P\right)
(D^c)+ \EE_{\nu_n}\left( \II_D f^2 _{M,n}\right) \leq c(\varphi)^2
(2\e-\e^2)+ \EE_{\nu_n} \left(f^2_{M,n}\right).
$$
By the same computations  as in  \cite{N}[Lemma 12], one gets
$$
\EE_{\nu_n} (f_{M,n}^2) \leq \frac{1}{2n}\sum _{j\in \ZZ}
\frac{1}{n} \varphi ^2 \left( \frac{j}{n}\right) ,
$$
implying that $\varlimsup _{n\uparrow\infty} E_{\nu_n} \left( f_n
^2\right) \leq c(\varphi ) ^2 (2\e-\e^2)$. Since $\e$ is arbitrary,
we get the thesis.
\end{proof}

\section{The random walk on $\ZZ$  with jump rates $\{c_k\}_{k\in\ZZ}$}\label{secpietra}

Let us   recall how  one can  express a 1D nearest--neighbor  random
walk as space--time  change of a 1D Brownian motion (see \cite{S}
for a detailful and more general discussion).

Let $B$ be a  Brownian motion with $\EE(B^2 (t))=t$, defined on some
probability space $\left(\mathbb{W} ,\mathbb{F}, \PP\right)$ (note
that in \cite{S} $\EE(B^2 (t))=2t$, thus changing the final results
of some factor $2$). Denote by $L (t,y)$ the local time of $B$.
Then, $\PP$--almost surely,
\begin{equation}\label{localino}
\int _a ^b L(t,y) dy = \int _0 ^t \II_{\{ a\leq B(s)\leq b\}} ds\,
\qquad \forall t\geq 0, \; \forall a\leq b \,.
\end{equation}

Let $\nu$ be a Radon measure on $\RR$ (i.e. $\nu$ is a Borel
positive measure on $\RR$, bounded on bounded intervals). We write
$\text{supp}(\nu)$ for the support of $\nu$ and we assume that
$\text{supp}(\nu)$ is unbounded from below and from above, namely
$$
  \inf \left(\text{supp}(\nu)\right)
=-\infty\,,\qquad
  \sup \left(\text{supp}(\nu)\right) =\infty\,.
  $$
For each $x\in\text{supp}(\nu)$ and $t\geq 0$, set
\begin{align}
& \psi (t|x,\nu )= \int _\RR L(t,y-x) \nu (dy) , \label{wwf1}\\
&  \psi ^{-1}(t|x,\nu )=\sup\left\{ s\geq 0\,:\, \psi (s|x,\nu )
\leq t \right\}.\label{wwf2}
\end{align}
Finally, we set
\begin{equation}\label{wwf3}
Z(t|x,\nu)= B \left( \psi ^{-1}(t|x,\nu )\right)+x.
\end{equation}

The process $\left(Z(t|x,\nu),\, t\geq 0 \right) $, defined on
$\left(\mathbb{W},\mathbb{F}, \PP\right)$ has paths in the Skohorod
space $D([0,\infty),\RR)$ endowed of the   Skohorod metric $d_S$.
Due to \cite{S}[Theorem 1] and \cite{S}[Corollary 1], the following
holds

\begin{prop}\label{pietra} \cite{S}
Let $\{\nu_n\}_{n\geq 0} $, $\nu$   be Radon measures on $\RR$ with
support unbounded from below and from above and  let  $x_n\in
\text{supp}(\nu_n) $ be a converging sequence with
$\lim_{n\uparrow\infty}x_n =x$. Suppose that: \\
  i) $\nu_n ([a,b])\to \nu ([a,b])$ for
all $a<b$ with
$\nu(\{a\})=\nu(\{b\})=0$,\\
 ii) if
$y_n \in \text{supp}(\nu_n) $ is a converging sequence as
$n\uparrow\infty$, then $\lim _{n\uparrow\infty}y_n\in
\text{supp}(\nu)$.

Then
\begin{equation}
 \lim_{n\uparrow\infty}
  d_S \left( Z (\cdot| x_n, \nu_n), Z(\cdot |x,\nu)\right)=0\,,\qquad
  \PP\text{ a.s. }
  \end{equation}
\end{prop}

Let us recall another consequence of the results in  \cite{S} (see
also \cite{KK}[Section 2]):
\begin{prop}\label{RW} \cite{S}
Let $\{x_k\}_{k\in \ZZ}$ satisfy  $\lim _{k\downarrow -\infty}
x_k=-\infty$, $ \lim _{k\uparrow \infty} x_k=\infty$. Fix positive
constants $\{w_k\}_{k\in \ZZ}$ and set $\nu =\sum _{k\in \ZZ } w_k
\d_{x_k}$. Then  $Z(\cdot|x_j,\nu)$ is the continuous--time random
walk on $\{x_k\}_{k\in \ZZ}$ starting in $x_j$ such that after
reaching site $x_k$ it remains in $x_k$ for an exponential time with
mean
$$
 2 w_k \frac{ (x_{k+1}-x_k)(x_k- x_{k-1})}{ x_{k+1}-x_{k-1}}
$$
and then it jumps to $x_{k-1}$, $x_{k+1}$ respectively with
probability
$$
 \frac{ x_{k+1}-x_k}{ x_{k+1}-x_{k-1} }\,\text{ and }\,
 \frac{ x_k-x_{k-1}}{ x_{k+1}-x_{k-1} }\,.
$$
\end{prop}
\subsection{The case $\a\in (0,\infty)$}
Recall the definition of  $S $ given in (\ref{sos}) and set  $
 S_n(x) = \frac{S(\lfloor xn\rfloor )}{n}$.  Due to condition (\ref{condi})
\begin{equation}\label{condi1}
 \lim _{k\uparrow \infty} \frac{S(\pm k)}{k} =\a,
 \end{equation}
 implying that
 $\lim_{n\uparrow \infty} S_n(x)=\a x $
 for all $x\in \RR$.

Set
$$
 \nu_n =\frac{1}{2n}\sum _{k\in \ZZ} \d _{S(k)/n}\,,\qquad \nu
 (dx) = \frac{1}{2\a} dx.
$$
 Due to Proposition
\ref{RW},
\begin{equation}\label{scalo}
 X_n (t|x) =\frac{1}{n} S ^{-1} \left( n Z(t| S_n(x),\nu_n ) \right) \,, \qquad x\in
 \ZZ/n\,,
\end{equation}
is the random walk on $\ZZ/n$ starting at $x$ with generator $ H_n
:\RR^{\ZZ/n}\to \RR^{\ZZ/n}$ where
\begin{equation}\label{maremma}
 H_n f \left(\frac{k}{n}\right)=n ^2 c_{k}\left(
 f\left(\frac{k+1}{n}\right)-f\left(\frac{k}{n}\right)\right)+
n^2  c_{k-1} \left(
 f\left(\frac{k-1}{n}\right)-f\left(\frac{k}{n}\right)\right).
\end{equation}

\begin{lemma}\label{semplice}
For all $a<b$
$$
\lim _{n\uparrow \infty } \nu_n ([a,b])= \nu ([a,b]).
$$
\end{lemma}
\begin{proof}
Since $\lim _{n\uparrow\infty}\nu_n (\{x\})=0$ for all $x\in \RR$,
it is enough to consider the cases where $a=0$ or $b=0$. We deal
with the former (the latter is similar).

 Trivially, $\nu_n
([0,b])= (\bar k +1) /2n$ where $\bar k = \max \{k\geq 0\,:\,
S(k)\leq n b\} $. Note that $\bar k =\bar k (n)$. Due to
(\ref{condi1}), given $\e>0$ there exists $k(\e)$ such that
\begin{equation}\label{maglia}
 \frac{k}{n} (\a-\e)\leq \frac{S(k)}{n}\leq \frac{k}{n} (\a+\e)
\qquad \forall n\geq 1, \forall  k\geq k(\e).
\end{equation}
Due to  (\ref{condi1}) $\lim _{n\uparrow\infty}\bar k (n)=\infty$,
therefore  we can assume the above expression to be  true for
$k=\bar k, \bar k +1 $. Since $ S(\bar k)/n\leq b$ and $S(\bar k
+1)/n>b$, (\ref{maglia}) implies
$$
\frac{b}{2(\a+\e)}-\frac{1}{2n}<\frac{\bar k }{2n}\leq
\frac{b}{2(\a-\e)}
$$
and therefore the thesis.
\end{proof}
 Due to  the above lemma, Proposition
\ref{pietra} holds for all sequences $x_n$, $n\geq 1$, such that
$x_n\in \ZZ/n$ and $x=\lim _{n\uparrow\infty}x_n$. Moreover, due to
(\ref{localino}),  $\psi (t|x, (2\a)^{-1} dx)= (2\a)^{-1} t $ thus
implying
$$Z(\cdot|
x, (2\a)^{-1} dx)=B(2\a t)+x \sim  \sqrt{2\a} B(t) +x ,
$$
where $X \sim Y $ means that the random variables $X,Y$ have the
same law.

The proof of the hydrodynamic limit will be based on the following
technical result:

\begin{prop}\label{vacanza} Suppose that $\a\in (0,\infty)$.
Fix $t>0$. Then for all $x\in \RR$,
\begin{equation}\label{mare}
X_n (t| x_n)\Rightarrow
  \sqrt{2/\a}  B(t)  +
x    \qquad \text{ as } n\to\infty,
\end{equation}
where $\Rightarrow $ denotes weak convergence and $x_n = \lfloor xn
\rfloor  /n$.
\end{prop}
\begin{proof}
By Proposition \ref{pietra}, since $S_n(x_n)\to \a x$,
$\PP$--almost surely
\begin{equation}\label{massa}
\lim _{n\uparrow\infty}
 d_{S}\left(\,
 Z(\cdot | S_n(x_n),\nu_n ) ,Z(\cdot| \a x, (2\a)^{-1} dx)\,
 \right)=0.
\end{equation}
Since $\PP$--almost surely  $Z\left(\cdot|\a x, (2\a)^{-1}
dx)\,\right)$ is  continuous, (\ref{massa}) implies that
$\PP$--almost surely
\begin{equation}\label{carrara}
\lim _{n\uparrow\infty }\sup_{0\leq s\leq T}
\left| Z(s| S_n(x_n),\nu_n ) -Z(s| \a x, (2\a)^{-1} dx)\right|
=0\,,\qquad \forall T>0.
\end{equation}
Fix $a\in \RR$. Since  $S $ is  increasing
 and due to (\ref{scalo}),
$$
\PP \left(  X_n (t|x_n)\leq a   \right)
=
\PP\left(  Z(t| S_n(x_n),\nu_n )     \leq S_n(a) \right).
$$
Due to (\ref{carrara}) with $T=t$ and since, given $\e>0$,  $\a
a-\e\leq S_n(a)\leq \a a+\e$  for $n$ large enough, we obtain that
 \begin{multline}
  \PP \left(Z (\cdot|\a x, (2\a)^{-1} dx)\leq \a a-\e\right) \leq
 \varliminf_{n\uparrow\infty}\PP \left(  X_n (t|x_n)\leq a
\right)\\
\leq \varlimsup_{n\uparrow\infty}\PP \left(  X_n (t|x_n)\leq a
\right)\leq \PP \left(Z (\cdot| \a x, (2\a)^{-1} dx)\leq  \a
a+\e\right).
\end{multline}
Due to arbitrariness of $\e$,
$$\lim _{n\uparrow \infty} \PP \left(  X_n (t|x_n)\leq a
\right)= \PP\left(Z (\cdot| \a x, (2\a)^{-1} dx)\leq  \a a
\right)=\PP(\sqrt{2/\a} B(t)+x\leq a).
$$
\end{proof}

\subsection{The case $\a=\infty$}
Fix $x\in \RR$ and define $S^{(n)}:\ZZ\to \RR$ as
$$
 S^{(n)} \left( \lfloor n x \rfloor +k\right)= S \left( \lfloor n x \rfloor +k\right)
-S \left( \lfloor nx\rfloor\right),\qquad \forall k\in \ZZ .
$$
Define  the measure $\nu_n$  as
$$
\nu _n = \frac{1}{2n} \sum _{k\in \ZZ} \d_{S^{(n)}\left( \lfloor n x
\rfloor +k\right)/n   }.
$$
Let $x_n= \lfloor n x \rfloor/n$. Since
$$
S^{(n)} \left( \lfloor n x \rfloor +k+1\right)- S^{(n)} \left(
\lfloor n x \rfloor +k\right)=\frac{1}{ c_{\lfloor n x \rfloor +k}},
$$
Proposition \ref{RW} implies that
\begin{equation}\label{annibale}
X_n(t|x_n)=\frac{ 1}{n}\left( S^{(n)}\right)^{-1} \left( n
Z(t|0,\nu_n) \right)
\end{equation}
 is the continuous--time random
walk on $\ZZ/n$ starting at $x_n$ with Markov generator $H_n$
defined in (\ref{maremma}).

\begin{prop}\label{supervacanza} Suppose that $\a=\infty$ and that  the assumptions
of Theorem \ref{lento} are  satisfied. Fix $t>0$. Then for all $x\in
\RR$,
\begin{equation}\label{mare1}
X_n (t| x_n)\Rightarrow  x    \qquad \text{ as } n\to\infty,
\end{equation}
where $\Rightarrow $ denotes weak convergence and $x_n = \lfloor xn
\rfloor  /n$.
\end{prop}
\begin{proof}
Given $n$ and $u<v<w$ in $\ZZ/n$, it is simple to build on a same
probability space    random walks $X'_n(\cdot|u)$, $X'_n(\cdot|v)$,
$X'_n(\cdot|w)$ having respectively the same law of
 $X_n(\cdot|u)$, $X_n(\cdot|v)$,  $X_n(\cdot |w)$ and such that
 \begin{equation}\label{coupling}
X'_n(s|u) \leq X'_n(s|v)\leq X'_n(s|w)\,\qquad \forall s\geq 0.
\end{equation}
To this aim consider a family of independent Poisson processes
$\left\{ N_k ^-(\cdot)\right\}_{k\in \ZZ}$, $\left\{ N_k
^+(\cdot)\right\}_{k\in \ZZ}$ such that $E\left(N_k ^-(t )\right)= E
\left(  N_k ^+(t)\right)= c_k n^2 t $
 for all $k\in \ZZ$.  The random walk on $\ZZ/n$ starting in a generic
point $x\in \ZZ/n$ can be realized as follows: arrived at a point
$k/n$ at time $t$ the particle waits until time  $s$ where
$$
s=\inf\left\{ u>t\,:\,      N_k ^-(t-) \not =  N_k ^-(t)  \text{ or }   N_k ^+(t-) \not =  N_k ^+(t) \right\}.$$
 At time $s$ the particle jumps  on the left if  $ N_k ^-(t-) \not =  N_k ^-(t) $ otherwise it jumps on the right.

Due to such a coupling it is enough to prove the thesis for $x\in
\QQ$. We first prove that for all $a>0$
\begin{equation}\label{urca1}
\lim _{n\uparrow \infty} \PP \left( X_n (t|x_n) > x_n + a \right)
=0.
\end{equation}
Due to
  (\ref{annibale}), it is enough to prove that
  \begin{equation}\label{urca2}
  \lim _{n\uparrow\infty} \PP \left( Z(t|0,\nu_n) > w_n \right) =
  0 ,
\end{equation}
where
$$
w_n = \frac{1}{n} S^{(n)} \left(
\lfloor n a \rfloor +\lfloor n x \rfloor \right).
$$
On $(\mathbb{W},\mathbb{F}, \PP )$  define the hitting time
$$
\t_y =\inf \left\{ s\geq 0 \, :\, B_s = y\right\}.
$$
Then the  reflection principle implies
\begin{equation}\label{kurt}
\PP (\t_y \leq s ) = 2\PP ( B_s\geq y ) = 1-\frac{1}{\sqrt{2\p}} \int _{ -y /\sqrt{s}}  ^{y /\sqrt{s}}
e^{-\frac{z^2}{2}} dz.
\end{equation}
Due to definition (\ref{wwf3}),  if $ Z(t| 0,\nu_n) > w_n $ then
$\psi^{-1} \left(t|0,\nu_n\right) > \t_{w_n} $, which implies that $
\psi \left(\t_{w_n}|0,\nu_n\right)\leq t $.  Since $\psi (\cdot
|0,\nu_n)$ is not decreasing, for all $\d>0$,
\begin{equation}\label{kurt1}
\PP \left( Z(t| 0,\nu_n) > w_n \right) \leq
\PP (\t_{w_n} <\d w_n ^2 ) + \PP \left( \psi (\d w_n^2  \,|\,
0,\nu_n)\leq t \right).
\end{equation}
The first addendum in the r.h.s. can be treated by means of   (\ref{kurt}):
\begin{equation}\label{vonne}
\lim _{\d\downarrow 0} \varlimsup _{n\uparrow \infty} \PP ( \t_{w_n}
< \d w_n ^2 ) =0 .
\end{equation}
Due to (\ref{localino}), the scaling property of Brownian motion and
since the local time is jointly continuous with probability $1$, one
gets for all $s\geq 0 $ that
\begin{equation}\label{annunzio}
L ( s, \cdot ) \sim L ( 1, \cdot / \sqrt{s}) \sqrt{s}.
\end{equation}
Since
$$
 \psi (\d w_n ^2 \,|\, 0,\nu_n) = \frac{1}{2n} \sum _{k\in \ZZ} L
 \left(
 \d w_n ^2 , \frac{ S^{(n)} (k) }{n} \right),
 $$
by taking $ s= \d w_n ^2 $ in (\ref{annunzio})  one gets
\begin{equation}\label{uovo}
\PP \left( \psi (\d w_n^2 \,|\, 0,\nu_n)  \leq t\right) = \PP ( Y_n
\leq t ),
\end{equation}
where
$$
Y_n =\frac{1}{2n } \sum _{k\in \ZZ} L \left (
1,
 \frac{ S^{(n)} (k) }{ \sqrt{\d } S ^{(n)} \left(
 \lfloor n a \rfloor + \lfloor n x \rfloor
 \right)
} \right) \frac{
 \sqrt{\d } S ^{(n)} \left(
 \lfloor n a \rfloor + \lfloor n x \rfloor
 \right)}{n}.
$$
Consider the event
$$
\Bb _{\rho ,\e} =\left \{ L(1,y)\geq \e\; \forall y\in
[0,\rho]\right\}.
$$
On $\Bb _{\rho,\e}$ it holds
$$
Y_n \geq \frac{\e \sqrt{\d} }{2n} c_n  ( \sqrt{\d} \rho  )
 \frac{
  S ^{(n)} \left(
 \lfloor n a \rfloor + \lfloor n x \rfloor
 \right)}{n}
,
$$
where
$$
c_n (\sqrt{\d}\rho ) = \left|\left\{ k\geq 0 \, :\, \frac{ S ^{(n)}
\left(
 k + \lfloor n x \rfloor
 \right)}{
 S ^{(n)} \left(
 \lfloor n a \rfloor + \lfloor n x \rfloor
 \right)} \leq \sqrt{\d}\rho
 \right\}\right|.
$$
Due to condition (H ) we can find a non negative  sequence $b_n$
such that
\begin{align}
&
 \lim _{n\uparrow \infty} \frac{b_n+1}{n}\,
\frac{ S ^{(n)} \left( \lfloor a n \rfloor +\lfloor x n \rfloor
\right) }{n}   =\infty , \label{forza2}\\
& \varlimsup _{n\uparrow \infty}
 \frac{
  S^{(n)} \left( b_n  + \lfloor n x \rfloor \right)
  }{
 S^{(n)} \left( \lfloor n a   \rfloor + \lfloor n x \rfloor \right)}
< \sqrt{\d} \rho, \label{forza1}
\end{align}
Due to (\ref{forza1}), $ c_n ( \sqrt{\d} \rho )\geq b _n +1$ for $n$
large enough. Therefore, on $\Bb _{\rho ,\e}$,
$$
Y_n \geq \frac{\e \sqrt{\d} (b _n+1) }{2n } \frac{
S ^{(n)} \left( \lfloor a n \rfloor +\lfloor x n \rfloor \right)
}{n}.
$$
The above estimate and (\ref{forza2}) imply that $\lim
_{n\uparrow\infty } Y_n =\infty$ on $\Bb _{\rho,\e}$. In particular,
\begin{equation}\label{friuli}
\varlimsup _{n\uparrow \infty}  \PP (Y_n\leq t) \leq \varlimsup
_{\rho \downarrow 0} \varlimsup _{\e\downarrow 0 } \PP \left( \Bb ^c
_{\rho ,\e} \right) .
\end{equation}
Since $L(1,0)>0 $ and $L(1,\cdot)$ is continuous with probability
$1$, the r.h.s. in the above expression is zero. This result
together with (\ref{kurt1}), (\ref{vonne}) and (\ref{uovo}) implies
 (\ref{urca1}). Similarly one can prove that
\begin{equation}\label{urca100}
\lim _{n\uparrow \infty} \PP \left( X_n (t|x_n) < x_n - a \right)
=0.
\end{equation}
(\ref{urca1}) and (\ref{urca100}) imply that $X_n (t
|x_n)\Rightarrow x $.

\medskip
\end{proof}

\begin{rem}\label{jolly}
Suppose  that
   there exists  an increasing sequence of positive integers $\{n_k\}_{k\geq 1}$  satisfying
   (\ref{tordo}),  (\ref{cca1}) and (\ref{cca2}). We claim that Proposition \ref{supervacanza} remains true.
    To this aim observe that the arguments of the above proof together with   (\ref{cca1}) and (\ref{cca2})   imply
    that
   \begin{equation}\label{superascoli}
   \lim _{k\uparrow \infty} \sup _{t\leq T } \PP \left ( \left| X_{n_k} (t| x_{n_k})-x\right|> a\right) =0,\qquad
   \forall T>0, a>0.
   \end{equation}
   For each integer $n$ set $ \underline n = n_K $, $\bar n = n_{K+1}$ where
   $   K=\sup \{ k\,:\, n_k\leq n \}.$
   Due to the coupling discussed   at the beginning of the above
   proof it holds
   $$
    X_{\underline  n}  (t (n/\underline  n )^2 ,x_{\underline n}
    )\, \underline n / n \leq X_n (t,x_n) \leq  X_{\bar n}  (t(n/\bar n )^2 ,x_{\bar n}
    )\, \bar n /n
    .
    $$
In particular, for
   each $a>0$,
   \begin{align}
   & \PP \left( X_n (t,x_n) > x+a  \right)\leq
    \PP \left( X_{\bar n}  (t(n/\bar n )^2 ,x_{\bar n} ) > (x+a) n /\bar n   \right), \label{asc124}\\
    &   \PP \left( X_n (t,x_n) < x-a  \right)\leq
    \PP \left( X_{\underline  n}  (t (n/\underline  n )^2 ,x_{\underline n} ) < (x-a) n /\underline n   \right).
       \label{asc125} \end{align}
   Note that due to (\ref{tordo}) $\lim _{n\uparrow \infty} n/\underline n =\lim _{n\uparrow \infty} n/\bar n=1$.  Therefore, (\ref{superascoli}), (\ref{asc124}) and (\ref{asc125}) imply (\ref{mare1}).
 \end{rem}


\section{Proof of Theorem \ref{Hyd} and Theorem \ref{lento} }\label{fine}
We  point out that in  both cases   $\a\in  (0,\infty)$ and $\a=\infty$
 it holds
 \begin{equation}\label{rumori}
 \frac{1}{n} \sum _{k\in \ZZ} \varphi\left(\frac{k}{n}\right)\sum
_{j\in \ZZ} p(t n^2, k,j) \h _j (0)=  \frac{1}{n} \sum _{j\in \ZZ}
\h _j (0) \EE\left(\varphi\left(X_n(t|j/n)\right)\right) ,
\end{equation}
where $X_n (t |\cdot )$ has been defined in (\ref{scalo}) for $\a\in
(0,\infty)$ and in (\ref{annibale}) for  $\a=\infty$.

 Let us first prove Theorem \ref{Hyd}.
  Let $\varphi \in C_c(\RR)$ and  set $g(x)= \EE\left(\varphi
\left( \sqrt{2/\a} B(t) +x \right)\right)$. Since $g\in C_0
(\RR)\cap L^1(\RR)$ and
$$
\int _\RR g(x)\rho _0 (x) dx = \int _\RR \varphi (x) \rho (x,t) dx,
$$
due to our assumption on $\mu_n$ we get
$$
\lim _{n\uparrow \infty} \mu_n \left( \left| \frac{1}{n} \sum _{k\in
\ZZ}  g\left(\frac{k}{n}\right) \h_k - \int _\RR \varphi(x) \rho
(x,t) dx \right|>\d \right)  =0.
$$
Due to the above limit,  Proposition \ref{gaetano} and
(\ref{rumori}),   in order to prove Theorem \ref{Hyd} it is enough
to show that
$$
\lim _{n\uparrow \infty}  \mu_n \left(  \frac{1}{n} \sum _{j\in
\ZZ}\left|\EE\left(\varphi\left(X_n(t|j/n)\right)\right)-
\EE\left(\varphi \left( \sqrt{2/\a} B(t) +j/n  \right)\right)\right|
\h _j (0) \right) =0
$$
Since the above expectation is bounded by
$$
\frac{1}{n} \sum _{j\in
\ZZ}\left|\EE\left(\varphi\left(X_n(t|j/n)\right)\right)-
\EE\left(\varphi \left( \sqrt{2/\a} B(t) +j/n
\right)\right)\right|,
$$
due to Scheff\'e Theorem (see the arguments in \cite{N} after
Statement 15 or the proof of Theorem \ref{lento} below) and the
uniform continuity of $\varphi$ it is enough to prove that
\begin{equation}\label{paganini}
\lim_{n\uparrow \infty}
\EE\left(\varphi\left(X_n (t|x_n )\right)\right)= \EE\left(\varphi
\left( \sqrt{2/\a} B(t) +x  \right)\right)\,,\qquad \forall x\in
\RR,
\end{equation}
where  $x_n=\lfloor x n \rfloor $. The above limit follows from
Proposition \ref{vacanza}, thus concluding the proof of Theorem
\ref{Hyd}.

\bigskip

Let us prove Theorem \ref{lento}. Due to (\ref{ricamo}), Proposition
\ref{gaetano} and  (\ref{rumori}), it is enough to prove that
$$
\lim _{n\uparrow \infty}  \mu_n  \left(  \left|
\frac{1}{n} \sum _{j\in \ZZ}
\EE\left(\varphi\left(X_n(t|j/n)\right)\right)\h _j-
\frac{1}{n} \sum _{j\in \ZZ} \varphi\left( \frac{j}{n}\right) \h _j
\right|
 \right) =0.
$$ Without loss of generality we
can assume $\varphi \geq 0$. Then  the above expectation is bounded
by
\begin{equation}\label{ranetta}
\int _{\RR} \left|f_n(x) - \varphi (x) \right| dx +\e_n
\end{equation}
where $\lim _{n\uparrow \infty} \e_n =0$  and
$$ f_n (x) =  \sum _{j\in \ZZ}
\EE\left(\varphi\left(X_n(t|j/n)\right)\right) \II_{\{x\in[ j/n,
(j+1)/n)\}}.
$$
In order to  conclude it is enough to apply the   same arguments
described in  \cite{N}   after Statement 15:
    $f_n\geq 0$, $\lim _{n\uparrow \infty} f_n (x) = \varphi
(x)$ for all $x\in \RR$ due to Proposition \ref{supervacanza} and
\begin{multline*}
\int _\RR f_n (x) dx = \frac{1}{n} \sum _{j\in \ZZ}
\EE\left(\varphi\left(X_n(t|j/n)\right)\right)= \frac{1}{n}  \sum
_{j\in \ZZ}\sum _{k\in \ZZ} p (tn^2, j, k)\varphi \left(
\frac{k}{n}\right) \\
= \frac{1}{n}  \sum _{j\in \ZZ}\sum _{k\in \ZZ} p (tn^2, k,
j)\varphi \left( \frac{k}{n}\right)= \frac{1}{n} \sum _{k\in \ZZ}
\varphi \left( \frac{k}{n}\right)\to \int_\RR \varphi (x)
dx.
 \end{multline*}
In particular we can apply Scheff\'e Theorem and get that the
integral in (\ref{ranetta}) goes to $0$, thus allowing to conclude
the proof of Theorem \ref{lento}.

\section{ Proof of Proposition \ref{lentone}}\label{isa1}

Since $c_k$ is in the domain of attraction of a $\nu$--stable law  we can write
(see \cite{EKM}[Theorem 2.2.8])
$$ \PP \left(  \frac{1}{c_k} > y \right) = L(y) y ^{-\nu},\qquad \forall y \geq 1,
$$
with $L$ slowly varying function. In particular (see
\cite{EKM}[Theorem A3.3]) $L$ can be written as
\begin{equation}\label{leo}
L(x) = h(x) \exp\left\{\int _1 ^x \frac{g(u)}{u} du \right\},\qquad
\forall x\geq 1,
\end{equation}
for suitable functions $h, g$ such that $h>0$, $\lim _{x\uparrow \infty} h(x) =
h_0\in (0,\infty)$, $\lim _{x\uparrow \infty} g(x) = 0 $.

The proof of Proposition \ref{lentone} is based on the following lemma:

\begin{lemma}\label{falo}
Let $Y_1, Y_2,\dots $ be a sequence  of i.i.d. random variables distributed as $1/c_0$ and
set
$$ M_n =\max\left\{ Y_1, Y_2, \dots , Y_n\right\},\qquad S_n =\sum _{k=1} ^n Y_k. $$
Then, given  $\b,\d,\g >0$,  there exist  positive constants $c_1,
c_2$ such that for all $n\geq 1$
\begin{align}
&  \PP\left(M_n \leq n^{\frac{1}{\nu}-\d} \right) \leq \exp \left\{ -c_1 n ^{\d\nu}\right\}, \label{fuoco1}\\
&   \PP \left( S_{\lfloor n^\b\rfloor }\geq n^\g\right) \leq c_2  n
^{\b-\nu \g+\d} . \label{fuoco2}
 \end{align}
\end{lemma}
\begin{proof}
Since $ 1- z\leq e ^{-z}$ for all $z\geq 0 $ and due to (\ref{leo})
$$
\PP\left(M_n \leq n^{\frac{1}{\nu}-\d} \right) \leq \left( 1 -
L\left( n^{\frac{1}{\nu}-\d} \right) n^{\d\nu -1 }\right) ^n \leq
\exp\left\{ -  L\left( n^{\frac{1}{\nu}-\d} \right) n^{\d\nu
}\right\}.
$$
Therefore (\ref{fuoco1}) holds with $c_1=\inf_{x\geq 1 } L(x)$ which is positive due to
representation (\ref{leo}).

In order to prove (\ref{fuoco2})   we point out that, due to  (\ref{leo}),
$$
\lim_{x\uparrow\infty}  L(x) / x^u=0, \qquad \forall u>0.
$$
Due to the above observation and
 since $ 1- z\geq   e ^{-2z}$ for $z$ small enough, we get for $n\geq n_0$
 \begin{multline}\label{anno}
\PP\left( M _{\lfloor n^\b\rfloor } \geq n ^\g \right)= 1 - \left(
1- L(n^\g) n^{-\nu \g}\right) ^{\lfloor n ^\b\rfloor  }\leq 1 - \exp
\left\{- 2  L(n^\g)  \lfloor n ^\b \rfloor n^{ -\nu\g}\right\} \\
\leq  2  L(n^\g)  \lfloor n ^\b \rfloor n^{ -\nu\g}
   \leq  n ^{\b -\nu\g+\d  } .
\end{multline}
Moreover, by integration by parts,
\begin{equation}\label{secolo}
 \EE \left (Y _1 \II _{ Y_1 \leq n ^\g } \right)\leq 1 + \int _1 ^{n^\g} x d \left ( L(x) x ^{-\nu}\right)\leq
 1+ L (n^\g) n^{\g(1-\nu)} \leq c_3  n^{\g(1-\nu)+\d}   .
 \end{equation}
In particular
\begin{equation}\label{nuovo}
\PP \left ( S_{\lfloor n ^\b\rfloor } \geq n ^\g ,\; M_{\lfloor
n^\b\rfloor } \leq n ^\g \right) \leq \PP \left( \sum
_{j=1}^{\lfloor n ^\b\rfloor } Y _j \II _{ Y_j \leq n ^\g } \geq n
^\g \right)\leq n ^{\b -\g} \EE \left (Y _1 \II _{ Y_1 \leq n ^\g }
\right)   \leq c_3    n ^{\b -\nu\g +\d}  .
\end{equation}
  Bounds (\ref{anno}) and (\ref{nuovo}) imply that
  $$
\PP  \left ( S_{\lfloor n ^\b\rfloor } \geq n ^\g\right) \leq \PP\left( M _{\lfloor n^\b\rfloor }
\geq n ^\g \right)+ \PP \left ( S_{\lfloor n ^\b\rfloor } \geq n ^\g ,\; M_{\lfloor  n^\b\rfloor }
\leq n ^\g \right)\leq c_2 n ^{\b -\nu\g+\d } .
$$

\end{proof}

In order to prove Proposition \ref{lentone} we distinguish between
the cases  $0<\nu <1/2$ and $1/2 \leq \nu <1$.

If $0<\nu <1/2$, then take $\d>0$ such that $2+\d<1/\nu$. Due to
(\ref{fuoco1}) and Borel--Cantelli lemma for almost all
environments, given $x,a\in \QQ$ with $a\not =0$, there exists a
positive constant $c$ such that
$$
S\left( \lfloor a n \rfloor +\lfloor x n \rfloor \right) - S\left(
 \lfloor x n \rfloor \right) > c n ^{ 1/\nu-\d},\qquad \forall n\geq 1 .
$$
The above estimate implies conditions (\ref{ca1}) and (\ref{ca2})
with $b_n =0$. Hence it is enough to apply Theorem \ref{lento}.

If $1/2 \leq \nu <1$, set
$$
\begin{cases}
\b= 2-1/\nu +2\d, \\
\g= 1/\nu -2\d/\nu,\\
\end{cases}
$$
with $\d>0 $  small enough to have $5\d <1/\nu -1$ and $2\d<1$. Then
$\b,\g>0$ and
$$
\b-\nu \g +\d = 1-1/\nu + 5\d <0.
$$
Set
 $n_k = \lfloor k^\rho \rfloor $ and $b_k =\lfloor n_k ^\b \rfloor
 $, where $\rho >0$ is large enough to have $ \rho (1-1/\nu +
 5\d)<-1$. Due to this choice and by (\ref{fuoco1}), (\ref{fuoco2})
 and Borel--Cantelli lemma for almost all environments the following
 holds:

 Given $x,a\in \QQ$ with $a \not =0$, there exist positive
 constants $c_1, c_2$  such that for all $k\geq 1$
 \begin{align}
& S\left( \lfloor a n_k \rfloor +\lfloor x n_k  \rfloor \right) -
S\left(
 \lfloor x n_k \rfloor \right) > c_1 n_k ^{1/\nu -\d} ,\\
& S\left( b_k  +\lfloor x n_k  \rfloor \right) - S\left(
 \lfloor x n_k \rfloor \right) <c_2 n_k ^\g.
 \end{align}
In particular, the l.h.s. of (\ref{cca1}) with fixed $k$ is bounded
from below by $ c_1 \lfloor n_k ^{\b}\rfloor n^{ 1/\nu -2 -\d }$
while the l.h.s. of (\ref{cca2}) with fixed $k$ is bounded from
above by $ (c_2/c_1) n_k ^{\g -1/\nu +\d }$. Due to our choice of
$\b,\g$ conditions (\ref{cca1}) and (\ref{cca2}) are satisfied.
Hence the thesis follows from the discussion after  Theorem
\ref{lento}

\appendix

\section{A counterexample}\label{aaa}

  As already
noted, (\ref{asso}) implies condition (\ref{condi}). We show here
that the inverse implication is false if $\a=\infty$.

 Consider the subsets $A, B,
C\subset \ZZ $ defined as
$$ A =
\cup _{n=0 }^\infty \left [ 2^{2n}, 2^{2n +1} \right)\cap \ZZ
,\qquad B = A\cup \left( -A\right) , \qquad C= \ZZ\setminus B.
$$
Define
$$ c_j =
\begin{cases}
1 & \text{ if } j\in B,\\
 e^{-|j|}
 & \text{ if } j \not \in B.
\end{cases}
$$
Then
$$
\lim _{n\uparrow \infty} \frac{1}{2^{2n}-1} \sum _{j= 2^{2n}} ^{
2^{2n+1}-2} \frac{1}{c_j} =1.
$$
In particular, $(\ref{asso})$ cannot hold  for $\a=\infty$, $x=1$,
$y=2$. Let us verify that (\ref{condi}) is fulfilled for
$\a=\infty$. Given $k> 1$,
\begin{equation}\label{casa}
\sum _{j=1 } ^{k-1} \frac{1}{c_j} \geq \sum _{j=1} ^{N_k} e^j
=\frac{ e^{N_k+1} -e }{e-1}
\end{equation}
where $N_k = |[1, k)\cap B^c |$. We claim that
\begin{equation}\label{stanca}
  \inf _{k> 1}  \frac{ N_k}{k}>0 .
\end{equation}
In fact, if $2 ^{2n}\leq k < 2^{2n +1} $ then
$$
[1,k) \cap B ^c = \cup _{u=1} ^n [ 2^{2u-1}, 2 ^{2u})\cap \ZZ,
$$
which implies that
 \begin{equation}\label{moltostanca}
  N_k = \sum _{u =1 } ^n 2^{2u-1} = \frac{ 2}{3}(2^{2n}-1).
 \end{equation}
For  $ 2^{2n +1}\leq k < 2^{2n+2}$ the above identity implies  that
 \begin{equation}\label{stanchissima}
 \frac{ 2}{3}(2^{2n}-1)   =   N_{ 2^{2n}}\leq   N_k \leq
 N_{2^{2n+2}}= \frac{ 2}{3}(2^{2n+2}-1).
 \end{equation}
 (\ref{stanca}) follows from   (\ref{moltostanca}) and (\ref{stanchissima}).

Due to (\ref{casa}) and (\ref{stanca}) we get
$$
\liminf _{k\uparrow \infty} \frac{1}{k} \sum _{j=1 } ^{k-1}
\frac{1}{c_j} \geq  \left(\inf  _{n>1} \frac{ N_n}{n}\right)
 \liminf _{N \uparrow \infty} \frac{ e^{N +1} -e
}{N(e-1)}=\infty.
$$
By symmetry one gets
$$
\lim _{k\uparrow\infty}\frac{1}{k}\sum _{j=1}^{k}
\frac{1}{c_{-j}}=\infty.
$$
In conclusion, $\{c_k\,:\, k\in \ZZ\}$ satisfies (\ref{condi}) with
$\a=\infty$ but it does not satisfy (\ref{asso}) with $\a=\infty$.

\bigskip
\section{Exclusion processes with bond disorder on
$\ZZ^d$}\label{alto}
 As already noted in the Introduction, the
results described in Section \ref{extnagy}  can be easily
generalized to higher dimension and by the arguments of Section
\ref{fine} one can prove the hydrodynamic limit of the exclusion
process with bond disorder on $\ZZ^d$  when having suitable information on the associated random
walk. In order to clarify this point, we  recall in this Appendix
the main  steps of our method leading to the proof of Theorem \ref{Hyd},
 which can be easily extended to higher
dimension.

\bigskip

To this aim we denote by $\Ee$ the set of non oriented bonds of
$\ZZ^d$ and consider the exclusion process on $\ZZ ^d$ with
generator
\begin{equation}
\Ll f (\eta)=\sum _{b\in \Ee} c_b \left( f(\eta^b)-f(\eta)\right),
\end{equation}
where   $c_b\geq 0$ and  $\eta ^b$ is the configuration obtained
from $\eta$ by switching the  values at the vertexes of $b$. Due to
a standard percolation argument applied to  the graphical
representation of exclusion processes, if $\sup _{b\in \Ee} c_b
<\infty$ then the above exclusion process is well defined.

It is simple to check that all the arguments and the
 results of Section \ref{extnagy}
  remain valid in higher dimension when suitably changing
the notation. In particular, the analogous of Proposition
\ref{gaetano} holds:  Given  $\d>0$, $t>0$, $\varphi \in C_c
(\RR^d)$ and given a sequence of probability measures $\mu_n$ on
$\{0,1\}^{\ZZ^d}$,
\begin{equation}\label{limitonebis}
\lim _{n\uparrow\infty} \PP _{\mu_n} \left( \left| \frac{1}{n^d}
\sum _{k\in {\ZZ}^d } \varphi\left(\frac{k}{n}\right) \h _k (tn^2)-
\frac{1}{n^d } \sum _{j\in {\ZZ}^d}  \h _j (0) P_t^n
\varphi\left(\frac{j}{n}\right) \right|>\d \right)=0,
\end{equation}
where $P_t^n\varphi(j/n) =\EE \left( \varphi \left( X_n
(t|j/n)\right)\right)$ and $X_n (t|j/n)$ is the random walk on
$\ZZ^d/n$ starting at $j/n$ with generator $H_n$ defined as the
$d$--dimensional version of (\ref{maremma}), namely
\begin{equation}\label{eq1974}
H_n f\left(j/n\right) =n^2 \sum _{k\in \ZZ^d\,:\, \|k-j\|_\infty =1
} c_{\{k,j\}}\left(f\left(k/n\right)-f\left(j/n\right)\right),\qquad
\forall j\in \ZZ ^d.
\end{equation}

 Suppose that
\begin{equation}\label{proteine}
\lim _{n\uparrow \infty} \mu_n\left(
\frac{1}{n^d} \left|
 \sum _{j\in \ZZ ^d}\eta_j(0)
P_t^n \varphi \left(\frac{j}{n} \right)-  \sum _{j\in \ZZ ^d}
 \eta_j ( 0) P_t\varphi \left(\frac{j}{n}\right) \right|\geq \d
\right) =0, \qquad \forall
\varphi \in C_c (\RR^d),
 \end{equation} where $P_t\varphi (x)=
\EE\left( \varphi(x+W(t))\right)$ and $W$ is a Brownian motion on
$\RR^d$ starting at $0$ (one could even consider other limiting
processes).
Trivially, (\ref{proteine}) holds if
\begin{equation}\label{vitamine}
\lim _{n\uparrow \infty}
\frac{1}{n^d} \sum _{j\in \ZZ ^d} \left|
P_t^n \varphi (j/n)- P_t\varphi (j/n)\right|  =0, \qquad \forall
\varphi \in C_c (\RR^d).
 \end{equation}
Note that in Section \ref{fine} we have derived
(\ref{vitamine})  from Proposition \ref{vacanza}  and
Scheff\'e Theorem (the same method works in all dimensions).

\bigskip
At this point, when
having  (\ref{proteine}) or even (\ref{vitamine}),
 the $d$--dimensional version
of Theorem \ref{Hyd} is easily proven   and
the hydrodynamic equation coincides with the linear
heat equation associated to $W$.

\bigskip

\noindent
{\bf Acknowledgements}. We thank D. Gabrielli and C. Landim  for useful discussions.

\bigskip

\noindent
{\bf Note added in proof.}
After completing this work,  we were aware of \cite{JL} where Theorem \ref{Hyd}
is proven by a different approach.

\end{document}